\documentclass{article}%\textheight 10.in
\newtheorem{Th}{Theorem}
\newtheorem{Cor}{Corollary}
\newtheorem{Prop}{Proposition}
\newtheorem{Lm}{Lemma}

\newtheorem{rem}{Remark}

\renewcommand{\P}{{\cal P}}

\newcommand{\R}{{\bf R}}
\newcommand{\N}{{\bf N}}

\newcommand{\Q}{{\bf Q}}
\newcommand{\ns}{\mbox{ns}}

\newcommand{\U}{{\cal U}}

\newcommand{\K}{{\cal K}}

\newcommand{\Ga}{\Gamma}

\newcommand{\e}{\varepsilon}

\renewcommand{\a}{\alpha}
\renewcommand{\b}{\beta}
\renewcommand{\c}{\gamma}
\renewcommand{\d}{\delta}
\newcommand{\s}{\sigma}

\newcommand{\la}{\langle}
\newcommand{\ra}{\rangle}
\newcommand{\st}{\mbox{st}}
\newcommand{\dom}{\mbox{dom}}
\newenvironment{proof}{{\par\noindent\bf Proof. }}{$\Box$}
\newcommand{\supp}{\mbox{supp}}
\renewcommand{\*}{\,^*\!}\renewcommand{\o}{\,^\circ\,\!\!}
\author{L.Yu.Glebsky, E.I.Gordon and C.J.Rubio}
\title{On approximation of topological groups by finite algebraic systems.~II}
\date{}
\usepackage{latexsym}
\begin{document}
%%%%%%%%%%%%%%%%%%%%%%%%%%%%%%%%%%%%%%%%%%%%%%%%%%%%%%%%%%%%%
\maketitle

\begin{abstract}

Recall that a locally compact group $G$ is called unimodular if the left Haar
measure on $G$ is
equal to the right one. It is proved in this paper that $G$ is unimodular
iff it is approximable
by finite quasigroups (Latin squares).

\end{abstract}

\section{Introduction}

This paper is a continuation of the paper \cite{GG}. We prove here the
statement that was formulated in \cite{GG}
as a conjecture, namely, the following theorem.

\begin{Th} \label{Main_th}
A locally compact group $G$ is unimodular iff it is approximable by finite
quasigroups
(see Definitions 2 and 3 in \cite{GG}).
\end{Th}

Notice that this theorem gives a solution of an old problem, formulated in
\cite{HR}: to  characterize the class of all
unimodular locally compact groups.

The sufficiency was proved in \cite{GG} (Corollary 1 of Proposition 5).
It was proved also that
any discrete group is approximable by finite quasigroups (Proposition 4). So in
this paper we prove only the following

\begin{Prop} \label{main_prop}
Any non-discrete locally compact unimodular group $G$ is approximable by
finite quasigroups.
\end{Prop}

The proof of this proposition that will be discussed in the paper
is rather complicated. We start with the case of a compact group
$G$ to outline the main ideas of this proof. The case of locally
compact groups requires some technical modification that will be
discussed at the end of the paper. A nontrivial combinatorics of
latin squares based on a generalization of one result of
A.J.W.Hilton \cite{Hil,Hil2} is involved in the proof. This
combinatorics is discussed in \cite{GC}.
%in the appendix.
As in \cite{GG} we also
use the language of nonstandard analysis in some proofs of this
paper. It allows to simplify essentially the proofs of Theorems~\ref{Partition}
and \ref{lpartition}. All theorems are formulated in standard language and we
hope that the main ideas of the proofs are understandable for the
readers non-familiar with nonstandard analysis.

\section{A proof of Proposition \ref{main_prop} for a compact group $G$.}

We assume in this section that $G$ is a non-discrete compact group.
All subsets of $G$ we deal with
are assumed to be measurable with respect to the Haar measure $\nu$
that is assumed to be normalized --- $\nu(G)=1$.
Recall that any compact group is unimodular and thus $\nu$ is left and
right invariant. Let $U$ be a neighborhood of the
unit in $G$ and $\P$ be a finite partition of $G$. We say that $\P$ is
$U$-fine if
$\forall P\in\P\ \exists g\in G\ (P\subseteq gU)$ and $\P$ is equisize
if all sets in $\P$ have the same Haar measure.
The following theorem will be proved in the next section.

\begin{Th} \label{Partition}
There exists a $U$-fine equisize partition for any neighborhood of
the unit $U\subseteq G$.
\end{Th}

Let $\P=\{P_1,\dots,P_n\}$ be a partition that satisfies the
conditions of this theorem for some $U$. Consider the following
three-index matrix $w=\la w_{ijk}\ |\ 1\leq i,j,k\leq n\ra$, where
$$
w_{ijk}=\mathop{\int\!\int}\limits_{G\times
G}\chi_i(xy^{-1})\chi_j(y)\chi_k(x)d\nu(x)d\nu(y), \eqno (1)
$$
here $\chi_m(x)=\chi_{P_m}(x)$ is the characteristic function of a set
$P_m$, $m\leq n$.

Obviously $w_{ijk}\geq 0$. Let $S=\supp\ w=\{\la i,j,k\ra\ |\ w_{ijk}>0\}$.

\begin{Lm} \label{Prop_w}
The three-index matrix $w_{ijk}$ has the following properties.
\begin{enumerate}
\item $\sum\limits_iw_{ijk}=\sum\limits_jw_{ijk}=
\sum\limits_kw_{ijk}=\frac 1{n^2}$
\item $S\subseteq\{\la i,j,k\ra\ |\nu(\ P_i\cdot P_j\cap P_k)>0\}$
\end{enumerate}
\end{Lm}

\begin{proof}
The statement 2 follows immediately from the Fubini's theorem.
To prove the statement
1 we need the following identity $\sum\limits_m\chi_m(t)=1$ for any
$t\in G$, that follows from the fact that $\P$ is a
partition of $G$.

Now it is easy to see that, for example,
$$
\sum\limits_jw_{ijk}=
\int\limits_G\chi_k(x)d\nu(x)\int\limits_G\chi_i(xy^{-1})d\nu(y)=
\int\limits_G\chi_k(x)d\nu(x)\int\limits_G\chi_i(y)d\nu(y)=
\nu(P_k)\nu(P_i).
$$
We used here the left invariance of $\nu$ and the unimodularity of $G$
that implies
$\int\limits_Gf(y)d\nu(y)=\int\limits _Gf(y^{-1})d\nu(y)$.
Since the partition $\P$ is equisize we have
$\nu(P_i)=\nu(P_k)=\frac 1n$.
\end{proof}

To motivate the following consideration let us use an analogy with two-index matrices. Recall that an $n\times n$
matrix $B=\|p_{ij}\|$ is called is called bistochastic if $p_{ij}\geq 0$ and
$\sum\limits_{i=1}^np_{ij}=\sum\limits_{j=1}^np_{ij}=1$. According to a well-known G.Birkhoff's theorem
(cf., for example, \cite{Ryser}) in this case $B$ is a convex hull of permutations - the matrices
that consist of zeros and ones and contain the unique one in an each row and in an each column.
Then there exists a permutation $T$ such that $\supp\ T\subseteq\supp\ B$. Assume for a moment that the similar fact
holds for the three-indexed matrices. We say that a three-index matrix is three-stochastic if it is nonnegative and the
sum of elements in each line is equal to one. We call a line any set $L$ of triples of elements of $\{1,\dots,n\}$
such that in all triples in $L$ two indexes are fixed and the third run over $\{1,\dots,n\}$. Notice that if $w_{ijk}$
satisfies Lemma \ref{Prop_w} then $n^2w_{ijk}$ is three-stochastic. So we assume that the following statement is true.

(A). If $w_{ijk}$ satisfies Lemma \ref{Prop_w} then there exists a matrix $\d_{ijk}$ that consists of zeros and ones,
contains the unique one in each line and such that $\supp\ \d_{ijk}\subseteq\supp\
w_{ijk}$.

By the properties of $\d_{ijk}$ it is easy to see that $\supp\ \d_{ijk}$ is the graph of the
operation $\circ$ on $\{1,\dots,n\}$ such that $i\circ j=k$ iff $\c_{ijk}=1$. Let us denote the
algebra $\{1,\dots,n\}$ with the operation $\circ$ by $Q$. Since for any $i$ and $k$ there exists the
unique $j$ such $\c_{ijk}=1$ and for any $j$ and $k$ there exists the unique $i$ that $\c_{ijk}=1$
we have that the left and right cancellation laws hold in $Q$ and thus $Q$ is a quasigroup.

Fix an arbitrary injection $\a:Q\to G$ such that for any $i\leq n$ holds $\a(i)\in P_i$.
Notice that if $i\circ j=k$ then $\la i,j,k\ra\in\supp\ \c_{ijk}\subseteq\supp\ w_{ijk}$ and thus
$P_i\cdot P_j\cap P_k\neq\emptyset$ by Lemma \ref{Prop_w} (2).

So we proved under assumption (A) the following
\begin{Lm} \label{injection}
For any neighborhood of the unit $U$ of a compact group $G$ and for any $U$-fine equisize partition
$\P$ of $G$ there exists a finite quasigroup $Q$ and an injection $\a:Q\to G$ such that
\begin{enumerate}
\item $\forall P\in\P\exists q\in Q\ (\a(q)\in P);$
\item $\forall q_1,q_2\in Q\ (\a(q_1)\in P_1\in\P\land
      \a(q_2)\in P_2\in\P\land \a(q_1\circ q_2)\in P_3\in\P
\Longrightarrow P_1\cdot P_2\cap P_3\neq\emptyset)$.
\end{enumerate}
\end{Lm}

Unfortunately the statement (A) is not true in general (see, for example
\cite{GC}
%%% Appendix
) and a proof of Lemma \ref{injection} is more difficult.
We will discuss it later in this section. At first let us show
how Lemma \ref{injection} implies Proposition~\ref{main_prop} for a
compact group $G$.

We will use the nonstandard criterion of approximability by finite
quasigroups - the Proposition 9
of \cite{GG}. According to this Proposition and using the fact that
any element of the nonstandard extension $\*G$ of a
compact group $G$ is nearstandard \cite{GG}, Proposition 8(3),
we have to prove the existence of a hyperfinite quasigroup
$Q$ and an internal map $\a:Q\to\* G$ that satisfy the following
conditions
\begin{itemize}
\item[(i).] $\forall g\in G\exists q\in Q\ (\a(q)\approx G);$
\item[(ii).] $\forall q_1,q_2\in Q\ (\a(q_1\circ q_2)\approx
\a(q_1)\cdot\a(q_2)).$
\end{itemize}

Fix an infinitesimal neighborhood of unit $U\subseteq \*G$ i.e.
$U$ is an internal open set and
$e\in U\subset\mu(e)$, where $\mu(e)$ is the monad of $\e$.
We may assume that $U$ is symmetric ($U=U^{-1}$) without loss
of generality.

By Theorem \ref{Partition} and the transfer principle there exist
a hyperfinite $U$-fine equisize partition $\P$ and thus there
exists a hyperfinite quasigroup $Q$ and an internal map $\a:Q\to
G$ that satisfy Lemma \ref{injection}. Let us show that $\la
Q,\a\ra$ satisfies the conditions (i) and (ii).

It is easy to see that all elements of any $X\in\P$ are
infinitesimally close to each other since $U$ is infinitesimal and
$\P$ is $U$-fine. Thus if $X,Y\in\P$ then all elements of $X\cdot
Y$ are infinitesimally close to each other.

Let $g\in G$. Since $\P$ is a partition of $\*G\supset G$ there exists
$P\in\P$ such that $g\in P$
By the condition (1) of Lemma~\ref{injection} there exist $q\in Q$
such that $\a(q)\in P$. So
$\a(q)\approx g$ and (i) is proved.

Let $\a(q_1)\in P_1,\ \a(q_2)\in P_2,\ \a(q_1\circ q_2)\in P_3$.
By the condition (2) of Lemma
\ref{injection} there exists $h\in P_1\cdot P_2\cap P_3$.
Then $\a(q_1)\cdot\a(q_2)\approx h$ and
$\a(q_1\circ q_2)\approx h$. This proves (ii) and
Proposition~\ref{main_prop}
for a compact group $G$.

We have only to prove Lemma \ref{injection}.

Let $Q$ be a quasigroup and $\s$ an equivalence relation on $Q$
which we will identify with a partition of $Q$ by equivalence
classes. So $\s=\{Q_1,\dots, Q_n\}$. Denote by $Q/\s$  the subset
of $\{1,\dots,n\}^3$ such that $\la ijk\ra\in Q/\s$ iff there
exist $q\in Q_i$ and $q'\in Q_j$ with $q\circ q'\in Q_j$. Notice
that if $\s$ is a congruence relation on $Q$ (i.e. it preserves
the operation $\circ$) then the introduced set is exactly the
graph of the operation in the quotient quasigroup $Q$ by $\s$ and
so we will call the set $Q/\s$ - a generalized quotient quasigroup
(gqq).

The following weakening of the statement (A) holds.

\begin{Th} \label{weak-birkhoff} Let a non-negative three indexes matrix
$u=\la u_{ijk}\ |\ 1\leq i,j,k\leq n\ra$ satisfy the following
condition
$$
\sum\limits_iu_{ijk}=\sum\limits_ju_{ijk}=\sum\limits_ku_{ijk}=l
$$
for some positive $l$. Then there exist a finite quasigroup $Q$
and its partition $\s=\{Q_1,\dots, Q_n\}$ such that the gqq
$Q/\s\subseteq\supp\ u$.
\end{Th}
This theorem easily follows from results of A.J.W.Hilton, \cite{Hil}.
See also \cite{GC} for a proof.
%See Appendix for the proof of Theorem \ref{weak-birkhoff}.

Now we are able to complete the proof of Lemma~\ref{injection}.

Let $\P=\{P_1,\dots,P_n\}$ be a $U$-fine equisize partition of
$G$. Consider the three indexes matrix $w$ defined by formula (1).
By Lemma~\ref{Prop_w} (1) the matrix $w$ satisfies the conditions
of Theorem~\ref{weak-birkhoff} and thus there exists a finite
quasigroup $Q$ and its partition $\s=\{Q_1,\dots, Q_n\}$ such that
$Q/\s\subseteq\supp w$. Consider an arbitrary injection $\a:Q\to
G$ such that $\a(Q_i)\subset P_i$ for any $1\leq i\leq n$. Since
all $P_i$ are infinite (our group $G$ is non-discrete) such
an injection $\a$ exists. Obviously the first condition of
Lemma~\ref{injection} holds. If $\a(q_1)\in P_1,\ \a(q_2)\in P_2,\
\a(q_1\circ q_2)\in P_3$ then $q_1\in Q_1,\ q_2\in Q_2$ and
$q_1\circ q_2\in Q_3$. Then $\la 123\ra\in Q/\s$ by the definition
of $Q/\s$ and thus $\la 123\ra\in\supp w$ by Theorem~\ref{weak-birkhoff}.
So $P_1\cdot P_2\cap P_3\neq\emptyset$ by
Lemma \ref{Prop_w}. $\Box$

\bigskip

\section{Proof of Theorem \ref{Partition}}

\bigskip

In this section again $G$ is a non-discrete compact group and $\nu$ the Haar measure on $G$ such that $\nu(G)=1$

We will use the following version of Marriage Lemma for finite
non-atomic measurable spaces due to Rado \cite{Rado}

\begin{Th} \label{Rado}
Let $S$ be a measurable space with a finite non-atomic measure
$\mu$, $\{S_1,\dots,S_n\}$ - a collection of subsets of $S$
such that $\bigcup\limits_{i=1}^nS_i=S$, $\la\e_1,\dots,\e_n\ra\in\R^n$,
$\e_i>0,\ \sum\limits_{i=1}^n\e_i=\mu(S)$.
Then the following two statements are equivalent:
\begin{enumerate}
\item there exists a partition $\{P_1,\dots, P_n\}$ of $S$ such that
$\mu(P_i)=\e_i,\
P_i\subseteq S_i,\ i=1,\dots,n$;
\item for any $I\subseteq\{1,\dots,n\}$ holds
$\mu\left(\bigcup\limits_{i\in I}S_i\right)\geq\sum\limits_{i\in I}\e_i$.
\end{enumerate}
\end{Th}

\begin{rem}
Indeed, Theorem \ref{Rado} is a very particular case of the theorem,
proved in \cite{Rado}.
\end{rem}

\begin{Lm} \label{HU}

For any neighborhood $U$ of the unit in $G$ there exists a finite set
$H\subset G$
such that
\begin{enumerate}
\item $HU=G$;
\item $\forall\ I\subseteq H\ \nu(IU)\geq\frac {|I|}{|H|}$.
\end{enumerate}
\end{Lm}

Theorem~\ref{Partition} follows immediately from Theorem
\ref{Rado} and Lemma \ref{HU}. Indeed, let $H$ satisfy Lemma
\ref{HU}, $H=\{h_1,\dots,h_n\}$. Consider the collection
$\{h_1U,\dots,h_nU\}$ of subsets of $G$ and put $\e_i=n^{-1},\
i=1,\dots,n$. Then the condition (2) of Theorem \ref{Rado} is
equivalent to the condition (2) of Lemma \ref{HU} and thus there
exists the partition $\P$ that satisfies the conditions of Theorem
\ref{Rado}. Obviously this partition satisfies Theorem \ref
{Partition}.

The remaining part of this section is dedicated to a proof of Lemma \ref{HU}.

Let $\U$ be the base of neighborhoods of the unit of $G$ that consists
of all symmetric
($V=V^{-1}$) neighborhoods. For $V\in\U$ put

$$
\Ga(V)=\bigcup\limits_{n=1}^{\infty}V^n.
$$

It is well known (cf. for example \cite{HR}) that $\Ga(V)$ is a
complete, and thus clopen, subgroup of $G$. Let
$\Ga_1(V)=\Ga(V),\Ga_2(V),\dots,\Ga_m(V)$ be a collection of all
left cosets of $\Ga(V)$. It is easy to see that for any $i\leq m$
and for any neighborhood of the unit $W\subseteq V$ holds:
$$
\Ga_i(V)W=\Ga_i(V). \eqno(2)
$$

\begin{Prop} \label{measure}
Let  $U\in\U$. Let $A\subseteq G$. If there exists coset $\Ga_i(U)$ such that $\Ga_i(U)\cap A\neq\emptyset$
and $\Ga_i(U)\backslash\overline{A}\neq\emptyset$, then
$$
\nu(\overline{A})<\nu(AU).    \eqno (3)
$$
\end{Prop}

\begin{proof}
 It is easy to see that $AU=\overline AU$ since $U$ is
open. Let us show that $AU\neq\overline A$. Indeed, if
$AU=\overline A$ then $AU^2=\overline AU=AU=\overline A$. By
induction $AU^n=\overline A$ and thus $A\Ga(U)=\overline A$. But
$\Ga_i(U)=g\Ga(U)$ for any $g\in \Ga_i(U)\cap A$ and thus
$\Ga_i(U)\subseteq\overline A$. The contradiction. So
$AU\setminus\overline A\neq\emptyset$ and since this set is open
we have $\nu(AU\setminus\overline A)>0$. This inequality implies
inequality (3) since $\overline A\subset AU$.
\end{proof}

In the remaining part of Lemma \ref{HU} we use the language of
nonstandard analysis (cf. \cite{GG}, where this language is
discussed).

We will need some modification of theorem about invariant
integral introduced in \cite{GG} (cf. formula (6) in \cite{GG}).

Let $L$ be a locally compact group, $C\subseteq L$ --- a compact set,
$U$ --- a relatively compact neighborhood of the unit.
Fix an arbitrary infinitesimal
neighborhood of the unit $O\subseteq\*L$ and an internal
$\*$compact set $K$ such that
$L\subseteq K\subseteq\*L$. Since $U$ and $C$ are standard it is easy
to see that the following
conditions hold:
$O\cdot O^{-1}\subseteq\*U;\ (\*C)^2\subseteq K;\ \*C\cdot\*U\subseteq K$.
It was shown in the proof of Lemma 1 of \cite{GG} that if under these
conditions a hyperfinite set
$H\subseteq\*L$ is an optimal left $O$-grid of $K$ then $H$ can be endowed
with a binary operation
$\circ$, making $H$ a left quasigroup that is a $\*(C,U)$-approximation of
$L$. Since this holds for all standard $C$ and $U$ we see that the left
quasigroup $H$ is a hyperfinite approximation of $G$.
Now applying Lemmas~2,3 and Proposition~5 of \cite{GG} we obtain
the following
\begin{Prop} \label{integral}
Let $L$ be a locally compact group, $V\subseteq L$ - a compact set with the
non-empty interior,
$O$ - an infinitesimal neighborhood of the unit in $\*L$,
$K$ - an internal compact set such that
$L\subseteq K\subseteq\*L$, $H\subseteq L$ - a hyperfinite set that is an
optimal left (right)
$O$-grid for $K$, $\Delta=|H\cap V|^{-1}$. Consider the functional
$\Lambda_V(f)$ defined by the formula
$$
\Lambda_V(f)=\o\left(\Delta\sum\limits_{h\in H}\*f(h)\right).
$$
Then $\Lambda_V(f)$  is a left (right) invariant finite positive functional
 on $C_0(L)$ and thus defines the left (right) Haar measure $\nu_V$ on $L$.
\end{Prop}

\begin{Cor} \label{boundary}
If under conditions of Proposition \ref{integral} $C$ is a compact and
$\nu_V(\partial C)=0$ then
$\nu_V(C)=\Lambda_V(\chi_C)$, in particular, if $\nu_V(\partial V)=0$ then
$\nu_V(V)=1$.
\end{Cor}

Corollary~\ref{boundary} follows immediately from Proposition~1.2.18 of
\cite{Gor} $\Box$.

\begin{Cor} \label{comp-integral}
Let $G$ be a compact group, $O$ - an infinitesimal neighborhood of the
unit in $\*G$ and $H\subseteq \*G$ an optimal left or right $O$-grid in
$\*G$ then the functional $\Lambda$ defined by the formula
$$
\Lambda(f)=\o\left(\frac 1{|H|}\sum\limits_{h\in H}\*f(h)\right).
$$
is a positive invariant finite functional on $C(G)$ and thus defines
the normalized Haar measure $\nu$ on $G$.
\end{Cor}

We are going to proof the existence of a hyperfinite
set $H\subset\*G$ such that for any standard neighborhood of the
unit $U\in\U$ holds
\begin{enumerate}
\item[(H1)] $H\*U=\*G$;
\item[(H2)] $\forall$ internal $I\subseteq H\ \*\nu(I\*U)\geq
\frac{|I|}{|H|}$.
\end{enumerate}

By the {\bf transfer principle} this implies Lemma \ref{HU}. Let
$O$ be an infinitesimal neighborhood of the unit in $\*G$ and $H$
- an optimal $O$-covering of $\*G$ (see \cite{GG}). We are going
to prove that $H$ satisfies the conditions (H1) and (H2).

Consider the  functional $\Lambda(f)$ defined in
Corollary~\ref{comp-integral}.
Then for any $f\in C(G)$ holds
$$
\Lambda(f)=\int\limits_Gfd\nu.
$$
Now it is easy to see that for any compact set $C\subseteq G$
holds
$$
\nu(C)\geq \Lambda(\chi_C),   \eqno (4)
$$
and for any open set $W\subseteq G$ holds
$$
\nu(W)\leq \Lambda(\chi_W).   \eqno (5)
$$

Indeed by the general theory of integration on locally compact
spaces (see for example \cite{Halmos})
$$
\nu(C)=\inf\{\Lambda(f)\ |\  f\in C(G),\ f\geq\chi_C\}.
$$
By the positivity of $\Lambda$ for any $f\geq\chi_C$ holds
$\Lambda(f)\geq \Lambda(\chi_C)$. This proves (4). To prove (5) it
is enough to apply (4) to the set $G\setminus W$.

{\bf Remark}.
The inequalities (4) and (5) hold for an arbitrary locally compact group
$G$ if an open
set $W$ is relatively compact. But inequality (5) in this case requires
a little bit more complicated
proof, using, for example, regularity of a Haar measure.

For an arbitrary (internal or external) set $I\subseteq G$ put
$\st(I)=\{\o x\ |\ x\in G\}$.

\begin{Prop} \label{Nonst}
Let $W\subseteq G$ be a standard open set, $U$ - a standard open
neighborhood of the unit in $G$ and $I\subseteq\*G$ - an internal
set. Then the following inclusions hold:
\begin{enumerate}
\item $I\*W\*U\supseteq \*(\st(I)W)$;
\item $\*(\st(I)U)\supseteq I$.
\end{enumerate}
\end{Prop}

\begin{proof}
1). We will prove the stronger inclusion
$\*(\overline{\st(I)W})\subseteq I\*W\*U$. We have
$$
\overline{\st(I)W}=\overline{\st(I)}\overline W=\st(I)\overline W,
\eqno(6)
$$
since $\st(I)$ is a closed set for any internal $I$ (this follows
from saturation - see, for example, \cite{Alb}). Let
$x\in\*(\overline{\st(I)W})$ then there exists
$b\in\overline{\st(I)W}$ such that $x\approx b$ since
$\overline{\st(I)W}$ is a compact set. By (6) the element $b$ can
be represented in the form $b=\o i\cdot a$ for some $i\in I$ and
$a\in\overline W$. By the nonstandard characteristic of the
closure of a standard set there exists $w\in\*W$ such that
$w\approx a$ and thus $x\approx iw$.
%By the definition of relation
%$\approx$
This implies that $x\in iw\*U\subseteq I\*W\*U$.

2) We will show that $I\subseteq \st(I)\*U\subseteq\*(\st(I)U)$.
Indeed, since $G$ is a compact set for any $i\in I$ there exists
$\o i\in\st (I)$ and since $i\approx\o i$, we have $i\in\o
i\cdot\*U\subseteq \st(I)\*U$.
\end{proof}

To complete the proof of Lemma~\ref{HU} we will prove that the
constructed hyperfinite set $H$ satisfies the conditions (H1) and
(H2). Let $U$ be an arbitrary (standard) neighborhood of the unit.
Since $HO=\* G$ and $O$ is an infinitesimal neighborhood of the
unit, and thus $O\subset\*U$ we see that the condition (H1) holds.
We have to prove only that $\*\nu(I\*U)\geq\frac{|I|}{|H|}$ for
any internal $I\subseteq H$.

Fix three neighborhoods of the unit $U_1,U_2$ and $U_3$ such that
$U_2\in\U$ and $U_1U_2U_3\subseteq U$. Put $A=\st(I)U_1$. There
are two possibilities.

1). The set $\overline A$ is the union of $k$ cosets of
$\Ga(U_2)$. Let $G=\bigcup\limits_{i=1}^m\Ga_i(U_2)$ be the
decomposition of $G$ on the left cosets of $\Ga(U_2)$ and
$\overline A=\bigcup\limits_{i=1}^k\Ga_i(U_2)$. By the equality
(2)
$$
\overline A\subseteq
AU_2\subseteq\left(\bigcup\limits_{i=1}^k\Ga_i(U_2)\right)\cdot
U_2=\bigcup\limits_{i=1}^k\Ga_i(U_2)=\overline A.
$$
So $AU_2=\overline A$ and thus $\nu(AU_2)=\frac k{m}$

Let $H_i=\*\Ga_i(U_2)\cap H$. Then all sets $H_i,\ i=0,\dots, m$
have the same cardinality. Indeed, since $O$ is infinitesimal and
thus $O\subseteq\*U_2$ we have $\*\Ga_i(U_2)O=\*\Ga_i(U_2)$ by (2).
So $H_iO\subseteq \*\Ga_i(U_2)$ and since $HO=G$ we have $H_iO=
\*\Ga_i(U_2)$. Now if $|H_i|<\frac 1{m}|H|$ for some $i$ then
using the left shifts of $H_i$ we obtain the $O$-cover of $G$ that
contains less number of elements than $H$, but this contradicts to
the optimality of $H$.

By Proposition~\ref{Nonst} (1) we have
$\*\nu(I\*U)\geq\nu(AU_2)=\frac k{m}$. By Proposition~\ref{Nonst} (2)
$I\subseteq \*A\cap H\subseteq\bigcup\limits_{i=1}^kH_i$.
So $|I|\leq\frac {k|H|}{m}$
and thus $\*\nu(I\*U)\geq\frac{|I|}{|H|}.$

2). The sets $A$ and $U_2$ satisfy Proposition \ref{measure}. Then
$\*\nu(I\*U)\geq \nu(AU_2)$ as above; $\nu(AU_2)>\nu(\overline A)$
by Proposition~\ref{measure}; $\nu(\overline A)\geq
\Lambda(\chi_A)=\o\left(\frac{|\*A\cap H|}{|H|}\right)$ by (4). But
$I\subseteq \*A\cap H$ as above. So, $\o(\*\nu(I\*U))>\o(\frac{|I|}{|H|})$
and, consequently, $\*\nu(I\*U)\geq\frac{|I|}{|H|}$
$\Box$

\bigskip

\section{Proof of Proposition~\ref{main_prop} in the general case}
\label{Sec_general_case}

\bigskip

In this section we consider the general case of a locally compact
unimodular group $G$. We denote by $\nu$ a (left and right
simultaneously) Haar measure and by $\U$ - the base of
neighborhoods of the unit that consists of all symmetric
relatively compact neighborhoods. Our proof of Proposition~\ref{main_prop}
is based on the following generalization of
Theorem~\ref{Partition}

\begin{Th} \label{lpartition}
For any neighborhood of the unit $U$ and any compact set $B\subset G$
there exist a compact set $C\supset B$ that has a $U$-fine equisize partition of $C$.
\end{Th}

We will prove this theorem in the next section.

Let $\P=\{P_1,\dots,P_n\}$ be a $U$-fine equisize partition of a
compact set $C$. Similarly to (1) consider the three indexes
matrix $w=\la w_{ijk}\ |\ 1\leq i,j,k\leq n\ra$, where
$$
w_{ijk}=\mathop{\int\!\int}\limits_{C\times
C}\chi_i(xy^{-1})\chi_j(y)\chi_k(x)d\nu(x)d\nu(y), \eqno (7)
$$

Let $S=\{\la i,j\ra\ |\ P_i\cdot P_j\subset C\}$.
\begin{Lm} \label{lProp_w} (Compare with Lemma \ref{Prop_w})
The three-index matrix $w_{ijk}$ has the following properties.
\begin{enumerate}
\item $\sum\limits_{i=1}^nw_{ijk},\sum\limits_{j=1}^nw_{ijk},
\sum\limits_{k=1}^nw_{ijk}\leq\frac{\nu(C)^2}{n^2}$
\item $\forall\ \la i,j\ra\in S\
\sum\limits_{k=1}^nw_{ijk}=\frac{\nu(C)^2}{n^2}$
\item $\forall \la i,j\ra\in S\exists k\ w_{ijk}>0$ and
$\forall \la i,j,k \ra w_{i,j,k}>0\Longrightarrow
                          \nu(\ P_i\cdot P_j\cap P_k)>0$.
\end{enumerate}
\end{Lm}

\begin{proof}
Notice that since $\P$ is a partition of $C$ we have
$\sum\limits_{a=1}^n\chi_a(t)=\chi_C(t)$ and since $\P$ is
equisize $\forall i\leq n\ \nu(P_i)=\frac{\nu(C)}{n}$

Now
$$
\sum\limits_{i=1}^nw_{ijk}=\mathop{\int\!\int}\limits_{C\times
C}\chi_C(xy^{-1})\chi_j(y)\chi_k(x)d\nu(x)d\nu(y)\leq
\mathop{\int\!\int}\limits_{G\times
G}\chi_j(x)\chi_k(y)d\nu(x)d\nu(y)=\frac{\nu(C)^2}{n^2}
$$
$$
\sum\limits_{j=1}^nw_{ijk}=\mathop{\int\!\int}\limits_{C\times
C}\chi_i(xy^{-1})\chi_C(y)\chi_k(x)d\nu(x)d\nu(y)\leq
$$
$$
\mathop{\int\!\int}\limits_{G\times
G}\chi_i(xy^{-1})\chi_k(x)d\nu(x)d\nu(y)=
\int\limits_G\chi_k(x)d\nu(x)\int\limits_G\chi_k(xy^{-1})d\nu(y)
$$
Due to the unimodularity of $G$ we have
$\int\limits_G\chi_k(xy^{-1})d\nu(y)=\int\limits_G\chi_k(y)d\nu(y)$
and thus
$\sum\limits_{j=1}^nw_{ijk}\leq \frac{\nu(C)^2}{n^2}$.
The third inequality in 1) can be proved similarly.

To prove the equality notice that since $P_i\cdot P_j\subseteq C$
the equality $\chi_i(xy^{-1})\chi_j(y)=1$ implies $\chi_C(x)=1$ and thus
$$
\sum\limits_{k=1}^nw_{ijk}=\mathop{\int\!\int}\limits_{C\times
C}\chi_i(xy^{-1})\chi_j(y)\chi_k(x)d\nu(x)d\nu(y)=
\mathop{\int\!\int}\limits_{C\times
C}\chi_i(xy^{-1})\chi_j(y)\chi_C(x)d\nu(x)d\nu(y)=
$$
$$
=\mathop{\int\!\int}\limits_{G\times
G}\chi_i(xy^{-1})\chi_j(y)d\nu(x)d\nu(y)=\frac{\nu(C)^2}{n^2}.
$$
In the last equality we used the right invariance of $\nu$.
The first part of statement 3) follows immediately from the statement
2) and the second - from
Fubini's Theorem.
\end{proof}

\begin{rem}
The proof of this Lemma is the only place in the proof of
Proposition~\ref{main_prop}, where the
assumption about the unimodularity of $G$ is used.
\end{rem}

Since the statements of Lemma~\ref{lProp_w} are weaker than those of
Lemma~\ref{Prop_w} we need a generalization of
combinatorial Theorem~\ref{weak-birkhoff}.

Let $\circ:\dom(\circ)\to Q$ be a partial binary operation on a set $Q$,
i.e. $\dom(\circ)\subseteq Q\times Q$.
We say that $Q$ is a partial quasigroup if for any $a,b\in Q$ each of
the equations $a\circ x=b$ and $x\circ a=b$ has no
more than one solution.

\begin{Lm} \label{Part-complete}
Any finite partial quasigroup $Q$ can be completed to a finite quasigroup,
 i.e. there exists a finite
quasigroup $(Q',\circ')$ such that $Q\subseteq Q'$ and
$\circ\subseteq\circ'$.
\end{Lm}

The proof of this lemma follows immediately from the fact that any
Latin subsquare can be completed to a Latin
square \cite{Ryser}. We used this fact in \cite{GG} to prove the
approximability
of discrete groups by finite quasigroups (Propsition 4 of \cite{GG}).

Let $\s$ is an equivalence relation on a partial quasigroup $Q$
that we identify with the partition $\{Q_1,\dots,Q_n\}$ of $Q$ by
$\s$-equivalence classes. Then the generalized quotient partial
quasigroup $Q/\s\subseteq\{1,\dots,n\}^3$ is defined exactly in
the same way as the generalized quotient quasigroup (see Section~2).

\begin{Th} \label{gen-birkhoff}
Let a non-negative three indexes matrix
$w=\la w_{ijk}\ |\ 1\leq i,j,k\leq n\ra$ and a set
$S\subset\{1,\dots,n\}^2$
satisfy the following conditions:
\begin{enumerate}
\item $\sum\limits_{i=1}^nw_{ijk},\sum\limits_{j=1}^nw_{ijk},
\sum\limits_{k=1}^nw_{ijk}\leq l$;
\item $\forall\ \la i,j\ra\in S\ \sum\limits_{k=1}^nw_{ijk}=l$
\end{enumerate}
for some positive real $l$.

Then there exists a finite partial quasigroup $(Q,\circ)$ and
a partition $\s=\{Q_1,\dots,Q_n\}$
of $Q$ that satisfy the following conditions:
\begin{enumerate}
\item the set
$\bigcup\limits_{\la i,j\ra\in S}Q_i\times Q_j\subseteq\dom(\circ)$;
\item the generalized partial quotient quasigroup $Q/\s\subseteq\supp\,w$
\end{enumerate}
\end{Th}

Theorem~\ref{gen-birkhoff} is proved in \cite{GC}.
% Appendix.
Now we are able
to complete the proof of Proposition~\ref{main_prop}. Similarly to
the case of compact groups we will use the nonstandard criterion
of approximability of locally compact groups by finite
quasigroups. We have to show that there exist a hyperfinite
quasigroup $(Q',\circ)$ and an internal injective map:$\a:Q'\to \*G$
such that
\begin{enumerate}
\item[i] $\forall g\in G\exists q\in Q'\ \a(q)\approx g$;
\item[ii] $\forall q_1,q_2\in\a^{-1}\left(\ns(\*G)\right)\
\a(q_1\circ q_2)\approx\a(q_1)\cdot\a(q_2)$.
\end{enumerate}

Recall that $\ns(\*G)$ is the set of all nearstandard elements of $\*G$.
Since $G$ is a locally compact group, there exists
an internal compact set $C\supseteq\ns(\*G)$. By Theorem~\ref{lpartition}
and the transfer principle we may assume that $C$ has a hyperfinite
$U$-fine equisize partition $\P=\{P_1,\dots,P_N\}$ for
some infinitesimal neighborhood $U$ of the unit in $\*G$.
Let $w=\la w_{ijk}\ |\ 1\leq i,j,k,\leq N \ra$ be an internal three indexes
matrix defined by formula (7) with this partition
$\P$. Notice that if $P_i$ contains at least one nearstandard point then
$P_i\subseteq\ns(\*G)$, and if $P_i,P_j\subseteq\ns(\*G)$ then
$P_i\cdot P_j\subseteq\ns(\*G)$. So if
$S\subseteq\{1,\dots,N\}^2$ is the set defined before Lemma~\ref{lProp_w} then
$$
S\supseteq\{\la i,j\ra\ |\ P_i,P_j\subseteq \ns\}.
$$

Let $(Q,\circ)$  be a hyprfinite partial
quasigroup and $\s=\{Q_1,\dots, Q_N\}$ --- its partition that satisfy the
conditions of Theorem~\ref{gen-birkhoff}. Let $Q'$ be a hyperfinite
quasigroup that completes $Q$
(see Lemma~\ref{Part-complete}).
Now consider an arbitrary internal injection
$\a:Q'\to\*G$ such that
$\a(Q'\setminus Q)\subseteq \*G\setminus C$, $\a(Q_i)\subseteq P_i$,
$i=1,\dots,N$.
The rest part of the proof is exactly the same as for the case of compact
group $\Box$.

\section{Proof of Theorem \ref{lpartition}}

Theorem \ref{lpartition} follows immediately from Theorem
\ref{Rado} and the following modification of Lemma \ref{HU}

\begin{Lm} \label{HUC}
For any compact set $B\subseteq G$ and for any neighborhood of the
unit $U$ there exist a compact set $C\supseteq B$ and a finite set
$F\subset C$ such that $C\subset FU$ and $\forall I\subseteq F\
\nu(IU\cap C)\geq\frac{|I|}{|F|}\nu(C)$
\end{Lm}

The proof of this lemma repeats mainly the proof of Lemma \ref{HU}
but requires some additional considerations.

Once again we assume without loss of generality that $U$ is
symmetric. This relatively compact symmetric neighborhood of the unit $U$ is
fixed
throughout this section. We say that a
set $S\subseteq G$ is $U$-disconnected if
there exists a set $A\subseteq S,\ A\neq\emptyset,\ A\neq S$ such
that $AU\cap S=A$. Otherwise $S$ is called $U$-connected.

\begin{Lm} \label{U-con}
Any set $K$, such that $U^n\subseteq
K\subseteq U^{n+1}$ for some $n\geq 0$, is $U$-connected.
\end{Lm}

\begin{proof}
Let $K$ satisfy conditions of the Lemma and be
$U$-disconnected.  So there exists a set $X\subset K$ such that
$\emptyset\neq X,\ XU\cap K=X,\ Y=K\setminus X\neq\emptyset$. Thus
$\{X,Y\}$ is a partition of $K$ and $XU\cap Y=\emptyset$. The last
equality implies that $X\cap YU^{-1}=\emptyset$ and thus $X\cap
YU=\emptyset$ since $U$ is symmetric. Thus $YU\cap K=Y$. Consider
the map $\Phi:2^K\to 2^K$ such that $\Phi(A)=AU\cap K$. This map
obviously has the following properties:
\begin{itemize}
\item if $A\subseteq B$ then $\Phi(A)\subseteq\Phi(B)$;
\item $\Phi^{n+1}(\{e\})=K$, where $e$ is the unit of $G$
\item $\Phi(X)=X,\ \Phi(Y)=Y$.
\end{itemize}

But $e\in X$ or $e\in Y$. This brings us to a contradiction.
\end{proof}

We say that a compact set $C$ is regular if it is equal to the
closure of its interior. The following lemma is a modification of
Proposition \ref{measure} for compact case.

\begin{Lm} \label{<}
Let $C$ be a regular compact $C_0$ - its interior. If $C_0$ is
$U$-connected, $A\cap C_0\neq\emptyset$, and $C\setminus\overline
A\neq\emptyset$ then $\nu(\overline A\cap C)<\nu(AU\cap C)$
\end{Lm}

\begin{proof}
Since $\overline A\subseteq AU$ holds $\nu(\overline
A\cap C)\leq\nu(AU\cap C)$. We have only to prove that this
inequality is strict. It is enough to prove that $\nu(\overline
A\cap C_0)<\nu(AU\cap C_0)$. It will be proved if we show that
$C_0\cap AU\neq C_0\cap \overline A$. Indeed in this case the set
$(C_0\cap AU)\setminus (C_0\cap\overline A)=C_0\cap
AU\setminus\overline A\neq\emptyset$. And since this set is open
holds $\nu((C_0\cap AU)\setminus(C_0\cap\overline A)>0$.

Since $U$ is open holds $AU=\overline AU$. We have
$$
C_0\cap\overline AU\supseteq C_0\cap(\overline A\cap
C_0)U\supseteq C_0\cap\overline A.
$$

Suppose that the last inclusion here is indeed the equality. Let
us show that $C_0\setminus (C_0\cap\overline
A)=C_0\setminus\overline A\neq\emptyset$. Indeed, if
$C_0\subseteq\overline A$ then $\overline C_0\subseteq\overline A$
and since $C$ is regular we have $C\subseteq\overline A$, but this
contradicts the condition $C\setminus\overline A\neq\emptyset$.
This implies that $C_0$ is $U$-disconnected (take $C_0\cap A$ for
$A$ and $C_0$ for $S$ in the definition of a $U$-disconnected
set). Contradiction.
\end{proof}

Proposition~\ref{Nonst} also needed to be modified for the case of
a locally compact group.

\begin{Lm} \label{cnonst}
Let $C\subseteq G$ be a compact set, $W\subseteq G$ - an open
relatively compact set, $U\subseteq G$ - a relatively compact
neighborhood of the unit and $I\subseteq \* C$ - an internal set.
Then $\*(\st(I)W\cap C)\subseteq I\*W\*U\cap\* C$ and
$I\subseteq\*(\st(I)U\cap C)$
\end{Lm}

The proof of this Lemma is exactly the same as the one of
Proposition~\ref{Nonst} $\Box$.

In the proof of Lemma \ref{HU} we used the functional $\Lambda(f)$
that is indeed the functional $I(f)$ introduced in \cite{GG} for a
locally compact group $G$. This functional is defined for all
bounded functions on $G$ and for $f\in C_0(G)$ holds
$I(f)=\int\limits_Gfd\nu$. For more general functions $f$ the last
equality may fail. It may fail even for the characteristic
function of a regular compact set $C$. But it follows from
Proposition~1.2.18 of \cite{Gor} that if $C$ is a regular compact
set and $\nu(\partial C)=0$ then $I(\chi_C)=\nu(C)$. In what
follows we need to deal with a compact set $C$ for which the last
equality holds. The following Proposition shows that we can
always find such a compact set big enough.

\begin{Prop} \label{m-zero}
For every $n\in\N$ there exists a regular compact set $K$ such
that $U^n\subseteq K\subseteq U^{n+1}$ and $\nu(\partial K)=0$.
\end{Prop}

\begin{proof}
Let $(K)_0$ denote the interior of a set $K$. Consider the family
$\K$ of all all compact sets $K$ such that
$U^n\subseteq K\subseteq U^{n+1}$. Consider the partial order
$\prec$ on $\K\times\K$ such that
$K_1\prec K_2$ iff $K_1\subseteq (K_2)_0$. Let $\Xi$ be a maximal chain
in $\K$ with respect to this partial order.
If $\Xi$ is uncountable then $Z$ contains at least one compact set
$K$ with $\nu(\partial K)=0$ since $\nu(U^{n+1})$ is
finite. So we may assume that $\Xi$ is countable.
There are three possibilities under this assumption:

1). There exist $X,Y\in\Xi$ such that $X\prec Y$ but there does
not exist $Z\in\Xi$ such that $X\prec Z\prec Y$. In this case due
to regularity of $G$ there exist an open set $W$ and a compact set
$K$ such that $X\subseteq W\subseteq K\subseteq (Y)_0$. Due to the
maximality of $\Xi$ either $K=X$ and thus $X$ is a clopen set, or
$K=Y$ and $Y$ is a clopen set. Since the boundary of a clopen set
is the empty set, we are done in this case.

2). The maximal chain $\Xi$ contains the maximal element $X$. Then
the similar consideration shows that $X$ is clopen.

3). The order type of $\Xi$ is either $\eta$ or $1+\eta$, where
$\eta$ is the order type of $\Q$. Let us show that this case is
impossible. Let us discuss the case of $\eta$, the case of
$1+\eta$ is absolutely similar.

Let $\Xi=\{X_{\a}\ |\ \a\in\Q\}$, and $X_{\a}\prec X_{\b}$ iff
$\a<\b$. Fix an arbitrary irrational number $a$ and put
$Y=\overline{\bigcup_{\a<a}X_{\a}}$. Then it is easy to see that
for all $\a<a$ one has $X_{\a}\prec Y$ and for $\a>a$ one has
$Y\prec X_{\a}$. This
contradicts the maximality of $\Xi$.

So we proved that there exists a compact set $K$ such that
$U^n\subseteq K\subseteq U^{n+1}$ and $\nu(\partial K)=0$. If this
$K$ is not regular we can consider $K'=\overline{(K)_0}$. It is well
known that $K'$ is always regular, $\partial K'\subseteq\partial
K$ and thus $\nu(\partial K')=0$ and $K'\subseteq U^{n+1}$. On the
other hand an open set $U^n\subseteq K$ and thus $U^n\subseteq
(K)_0\subseteq K'$
\end{proof}

The proof of Lemma~\ref{HU} used the finite decomposition of the
compact group $G$ by the cosets of subgroup
$\Gamma(U)=\bigcup\limits_{n=1}^{\infty}U^n$. In general case this
subgroup may not be compact, but it is complete and thus clopen.
Also the number of cosets may not be finite. So we need the
following modification of the mentioned decomposition.

\begin{Lm} \label{decomp}
For any symmetric relatively compact neighborhood of the unit $V$
and for any compact set $B\subseteq G$ there exist a regular compact set
$C'$ and a finite set $\{g_1,\dots, g_n\}\subset G$ such that
\begin{itemize}
\item the interior $(C')_0$ of $C'$ is $V$-connected;
\item $\nu(\partial C')=0$;
\item $B\subseteq\bigcup\limits_{i=1}^ng_iC'$;
\item if $i\neq j$ then $g_iC'V\cap g_jC'V=\emptyset$
\end{itemize}
\end{Lm}

\begin{proof}
Consider the decomposition of $G$ into the family (maybe infinite) of
left cosets of
$\Gamma(V)$. Since the cosets are clopen sets there exist only finitely
many cosets
$\Gamma_1=g_1\Gamma(V),\dots,\Gamma_n=g_n\Gamma(V)$ such that
$\Gamma_i\cap B\neq\emptyset,\ i=1,\dots,n$ and there exists an
$m\in\N$ such that for all $i\leq n$ holds
$\Gamma_i\cap B\subseteq g_iV^m$.
By Proposition \ref{m-zero} there exist regular compact set
$C'$ with $\nu(\partial C')=0$ such that
$V^m\subseteq C'\subseteq V^{m+1}$. By Lemma \ref{U-con}
the set $C'$ is $V$-connected. It is easy to see that $C'$ satisfies
all other conditions of this lemma.
For example, $g_i\Gamma(V)\cdot V= g_i\Gamma(V)$ and thus the last
condition holds.
\end{proof}

We are able now to complete the proof of Lemma~\ref{HUC}.
Similarly to the proof of Lemma~\ref{HU}
we will show that there exist a compact set $C\supseteq B$ and
a hyperfinite set
$F\subset \*C$ such that $\* C\subset F\* U$ and for any internal
$I\subseteq F\ \*\nu(I\*U\cap \* C)\geq\frac{|I|}{|F|}\nu(C)$.
Lemma~\ref{HUC} follows from this statement by Transfer
Principle.

Let $U=U_1U_2U_3$. We may assume that $U_2$ is symmetric without
the loss of generality. Let
$C'$ and $g_1,\dots,g_n$ satisfy conditions of Lemma \ref{decomp}
for $B$ and $V=U_2$. Put
$C=g_1C'\cup\dots\cup g_nC'$. Consider an internal compact $K$
satisfying $\Gamma (U_2)\subseteq K\subseteq\*\Gamma(U_2)$.

There exists a hyperfinite set $M$ such that
\begin{itemize}
\item $\{g_1,\dots,g_n\}\subseteq M$;
\item $G\subseteq MK$;
\item for any $m_1\neq m_2\in M$ holds $m_1K\*U_2\cap m_2K\*U_2=\emptyset$.
\end{itemize}

Indeed, let $G\subseteq X\subseteq \*G$ be an internal compact set.
Let $D=\{E\in\*G/\*\Gamma(U_2)\ |\ E\cap X\neq\emptyset\}$.
Then $D$ is hyperfinite and for any $g\in G$ one has
$\*(g\Gamma(U_2))\in D$. Consider
an internal set $M'$ of representatives of the internal family
$D\setminus\{g_1\*\Gamma(U_2),\dots,g_n\*\Gamma(U_2)\}$ and
put $M=M'\cup\{g_1,\dots,g_n\}$. If $m_1\neq m_2\in M$ then
$m_1\*\Gamma(U_2)\*U_2\cap m_2\*\Gamma(U_2)\*U_2=\emptyset$ since
$\Gamma(U_2)U_2=\Gamma(U_2)$ and thus $m_1K\*U_2\cap m_2K\*U_2=\emptyset$.

Let $O\subseteq\*G$ be an infinitesimal neighborhood of the unit and
a hyperfinite set $H$ be
an optimal $O$-grid for $K$. Then obviously $MH$ is an optimal $O$-grid
for $MK$.
Put $F=MH\cap C$. By Proposition~\ref{integral} the functional
$$
\Lambda(f)=\o\left(\frac 1{|F|}\sum\limits_{x\in MH}\*f(x)\right)
$$
restricted on $C_0(G)$ is an invariant functional, which induces the Haar
measure $\nu_C$ on $G$
such that $\nu_C(C)=1$, for it is easy to see that $\nu(\partial C)=0$.
In what follows we identify $\nu_C$ and $\nu$.

Fix an arbitrary internal $I\subseteq F$. We have to prove that
$\*\nu(I\cdot\*U\cap\*C)\geq\frac{|I|}{|F|}$.

Let $A=\st(I)\cdot U_1$. There are two possibilities.

1). The set $\overline{A}\cap C$ is the union of $k\leq n$ sets
$g_iC'$, say,
$\overline{A}\cap C=\bigcup\limits_{i=1}^kg_iC'$.
Then $A\cdot U_2\cap C=\overline{A}\cap C$.
Indeed
$$
\overline {A}\cap C=\overline{A\cap C}\subseteq
(A\cap C)\cdot U_2\cap C\subseteq
(\bigcup\limits_{i=1}^kg_iC'U_2)\cap C =\overline{A}\cap C.
$$

The last equality holds since $g_iC'U_2\cap g_jC'=\emptyset$ for $i\neq j$.

So we have $\nu(AU_2\cap C)=\frac kn.$

Let $F_i=F\cap g_i\*C'$. Then
$F_i=MH\cap g_i\*C'=g_iH\cap g_i\*C'=g_i(\*C'\cap H)$. Thus all $F_i$ have
the same cardinality and so $|F_i|=\frac {|F|}n$.
By Lemma~\ref{cnonst} $I\subseteq \st(I)\*U_1\cap \*C$, thus
$|I|\leq\frac{k|F|}n$.

Again by Lemma \ref{cnonst}
$\*(\st(I)\cdot U_1\cdot U_2\cap C)\subseteq I\*U\cap C$, thus
$$
\*\nu(I\*U\cap \*C)\geq\nu(\st(I)\cdot U_1\cdot U_2\cap C)=
\frac kn\geq\frac{|I|}{|F|}.
$$

2). For some $i\leq n$ the sets $\st(I)U_1$, $g_iC'$ and $U_2$ satisfy
the conditions of Lemma~\ref{<}
for $A,C$ and $U$ respectively. Using Lemma~\ref{<}, Lemma~\ref{cnonst} and
inequality (4) we obtain:
$$
\*\nu(I\cdot\*U\cap C)\geq\nu(\st(I)\cdot U_1\cdot U_2\cap C)>
\nu(\overline{\st(I)\cdot U_1\cap C})\geq
$$

$$
\geq\Lambda(\chi_{\st(I)\cdot U_1\cap C})=
\o\left(\frac{|(\st(I)\cdot U_1\cap C)\cap F|}{|F|}\right)\geq
\o\left(\frac{|I|}{|F|}\right)\approx\frac{|I|}{|F|}\
\Box
$$
\section{Approximation of unimodular groups by loops}

In this section we sketch a proof of a result a little bit
stronger than Theorem 1, namely, the following

\begin{Th}\label{loop}
Any locally compact abelian group $G$ is approximable by finite
loops.
\end{Th}

Recall that an element $e$ of a quasigroup $(Q,\circ)$ is called
{\it the unity} if $\forall a\in Q\, a\circ e=e\circ a=a$. A
quasigroup with the unity is called {\it a loop}.

Recently, Milo\v{s} Ziman \cite{Ziman} proved that any discrete
group is approximable by loops, so we have to prove only the
following

\begin{Prop} \label{prop-loop}
Any non-discrete locally compact unimodular group is approximable
by finite loops.
\end{Prop}

To prove this proposition we need the following modifications of
Lemma 4 and Theorem 6.

\begin{Lm} \label{loop1}
The three-index matrix $w_{ijk}$ has the following properties.
\begin{enumerate}
\item $\sum\limits_{i=1}^nw_{ijk},\sum\limits_{j=1}^nw_{ijk},
\sum\limits_{k=1}^nw_{ijk}\leq\frac{\nu(C)^2}{n^2}$
\item $\forall\ \la i,j\ra\in S\
\sum\limits_{k=1}^nw_{ijk}=\frac{\nu(C)^2}{n^2}$; \hfill
 $\forall\ \la i,k\ra\in S'\
\sum\limits_{j=1}^nw_{ijk}=\frac{\nu(C)^2}{n^2}$; \hfill
 $\forall\ \la j,k\ra\in S''\
\sum\limits_{i=1}^nw_{ijk}=\frac{\nu(C)^2}{n^2}$;
\item
$\forall \la i,j,k \ra\; w_{ijk}>0 \Longrightarrow\nu(\ P_i\cdot
P_j\cap P_k)>0$,
\end{enumerate}
where $S'=\{\la i,k\ra\; |\; P_i^{-1}\cdot P_k\subseteq C\}$,
$S''=\{\la j,k\ra\; |\; P_k\cdot P_j^{-1}\subseteq C\}$.
\end{Lm}

The proof is the same as the one of Lemma~\ref{lProp_w} (see
Section 4).

\begin{Th} \label{loop2}
Let a non-negative three indexes matrix $w=\la w_{ijk}\ |\ 1\leq
i,j,k\leq n\ra$ and sets $S,S',S''\subset\{1,\dots,n\}^2$ satisfy
the following conditions:
\begin{enumerate}
\item $\sum\limits_{i=1}^nw_{ijk},\sum\limits_{j=1}^nw_{ijk},
\sum\limits_{k=1}^nw_{ijk}\leq l$;
\item $\forall\ \la i,j\ra\in S\ \sum\limits_{k=1}^nw_{ijk}=l$;\hfill
 $\forall\ \la i,k\ra\in S'\ \sum\limits_{j=1}^nw_{ijk}=l$;\hfill
 $\forall\ \la j,k\ra\in S''\ \sum\limits_{i=1}^nw_{ijk}=l$.
\end{enumerate}
for some positive real $l$.

Then there exists a finite partial quasigroup $(Q,\circ)$ and a
partition $\s=\{Q_1,\dots,Q_n\}$ of $Q$ that satisfy the following
conditions:
\begin{enumerate}
\item the set
$\bigcup\limits_{\la i,j\ra\in S}Q_i\times
Q_j\subseteq\dom(\circ)$;
\item equation $a\circ x=b$ ($x\circ a=b$) has a solution for
$\la a,b\ra\in \bigcup\limits_{\la i,j\ra\in S'}Q_i\times Q_j$
(for $\la a,b\ra\in\bigcup\limits_{\la i,j\ra\in S''}Q_i\times
Q_j$);
\item the generalized partial quotient quasigroup $Q/\s\subseteq\supp\,w$
\end{enumerate}
\end{Th}
The proof is an easy modification of the one of
Theorem~\ref{gen-birkhoff} (see \cite{GC}).

Using Lemma \ref{loop1} and Theorem \ref{loop2} one immediately
obtains that the quasigroup $Q'$ and the map $\a:Q'\to\* G$
constructed in the proof of Proposition 1 for the general case
(see the very end of Section 4) satisfy the following condition.

\bigskip

{\it (I) If $\a (x),\a (z)\in\ns$ and ($x\cdot y =z$ or $y\cdot
x=z$) then $\a (y)\in\ns$}.

\bigskip

Now we introduce a new loop operation $*$ on $Q'$ such that
$(Q',*)$ approximate $G$ with the same $\a$.

Construction:
\begin{itemize}
\item Take $q_0\in Q'$, such that $\a (q_0)\approx e$ ($e\in G$ is the unity).
\item Find permutation $a:Q'\to Q'$ such that $q_0\circ a(x)=x$.
By the property (I) $a(x)\in\ns$ if and only if $x\in \ns$. So, if
$x\in\ns$ then $a(x)\approx x$.
\item Find permutation $b:Q'\to Q'$ such that
$b(x)\circ a(q_0)=x$. By the same argument $b(x)\approx x$ for
$x\in\ns$. It easy to check that $b(q_0)=q_0$.
\item Define operation $x*y=b(x)\circ a(y)$. It is easy to see that $(Q',*)$
is a loop with the unity $q_0$ and $(Q',*)$ with $\a$ approximates
$G$.
\end{itemize}

This proves Proposition \ref{prop-loop} $\Box$

\bigskip

Instituto de Investigaci\'on en Communicaci\'on Optica de Universidad Aut\'onoma de San
Luis Potos\'\i, M\'exico

Eastern Illinois University, USA

\bigskip

IICO-UASLP

Av. Karakorum 1470

Lomas 4ta Secci\'on

San Luis Potos\'{\i} SLP 7820

M\'exico

Phone: 52-444-825-0892 (ext. 120)

e-mail:glebsky@cactus.iico.uaslp.mx\\

\bigskip

Mathematics Department 1-00036

Eastern Illinois University

600 Lincoln Avenue

Charleston, IL 61920-3099

USA

Phone: 1-217-581-6282

e-mail: cfyig@eiu.edu

\bigskip

1991 {\it Mathematics Subject Classification}. Primary 26E35, 03H05; Secondary 28E05, 42A38

\end{document}